\documentclass{amsart}

\usepackage{amsfonts}
\usepackage{amssymb}
\usepackage{mathabx}
\usepackage{amscd}
\usepackage{pictexwd,dcpic}
\usepackage{graphicx}
\usepackage{hyperref}
\hypersetup{
    colorlinks=true,
    linkcolor=black,
    filecolor=magenta, 
    urlcolor=cyan,
    citecolor=black,}
\usepackage[nameinlink,capitalize]{cleveref}

\title[Fiber cones and defining ideal of Rees algebras of modules]{Cohen-Macaulay fiber cones and \\ defining ideal of Rees algebras of modules}
\author{Alessandra Costantini}
\address{Department of Mathematics, University of California Riverside, Riverside CA 92521, USA}
\email{alessanc@ucr.edu}
\date{}	

\newtheorem{thm}{Theorem}[section]

\newtheorem{tdefn}[thm]{Theorem and Definition}
\newtheorem{prop}[thm]{Proposition}
\newtheorem{set}[thm]{Setting}
\newtheorem{notat}[thm]{Notation}
\newtheorem{lemma}[thm]{Lemma}
\newtheorem{cor}[thm]{Corollary}
\newtheorem{rmk}[thm]{Remark}

\begin{document}

\maketitle 

  \begin{abstract}
     Generic Bourbaki ideals were introduced by Simis, Ulrich and Vasconcelos in \cite{SUV2003} to study the Cohen-Macaulay property of Rees algebras of modules. In this article we prove that the same technique can sometimes be used to investigate the Cohen-Macaulay property of fiber cones of modules and to study the defining ideal of Rees algebras. This is possible as long as the Rees algebra of a given module $E$ is a deformation of the Rees algebra of a generic Bourbaki ideal $I$ of $E$. Our main technical result provides a deformation condition that in fact extends the applicability of generic Bourbaki ideals to situations not covered in \cite{SUV2003}.
  \end{abstract}

\section{Introduction}
  
  In this paper we study the Rees algebra $\,\mathcal{R}(E)\,$ and the fiber cone $\,\mathcal{F}(E)\,$ of a finitely generated $R$-module $E$, where $R$ is either a Noetherian local ring or a standard graded ring. The Rees algebra $\,\mathcal{R}(E)\,$ is defined as the symmetric algebra $\,\mathcal{S}(E)\,$ modulo its $R$-torsion submodule. The fiber cone of $\,\mathcal{F}(E)\,$ is then obtained by tensoring the Rees algebra with the residue field $k$ of $R$. 

  Much of the motivation to study Rees algebras and fiber cones comes from algebraic geometry. Indeed, Rees algebras arise for instance as homogeneous coordinate rings of blow-ups of schemes along one or more subschemes, or as bihomogeneous coordinate rings of graphs of rational maps between varieties in projective spaces. Correspondingly, the fiber cone is the homogeneous coordinate ring of the special fiber of the blow-up at the unique closed point, or of the image of the given rational map. In many situations, these are Rees algebras and fiber cones of modules which are not ideals, for instance when considering a sequence of successive blow-ups of a scheme along two or more disjoint subschemes, or the Gauss map on an algebraic variety.
  In addition, Rees algebras and fiber cones of modules arise when studying how algebraic varieties and morphisms change as they vary in families, in connection with the theory of multiplicity and Whitney equisingularity \cite{{Teissier1},{Teissier2}}. Algebraically, this relates to the study of integral dependence of ideals and modules \cite{{NR},{ReesAmao},{Rees},{SUVmult},{KR},{KT},{FM},{UVj},{UVepsilon}}.

 In this paper we have two main goals. The first is to understand the Cohen-Macaulay property of the fiber cone $\,\mathcal{F}(E)\,$ of a module $E$. In \cite{SUV2003}, Simis, Ulrich and Vasconcelos introduced the notion of \emph{generic Bourbaki ideals} to reduce the study of the Cohen-Macaulay property of Rees algebras of modules to the case of ideals, exploiting the remarkable fact that torsion-free module of rank one are isomorphic to ideals of positive grade (see \cref{SecPrelim} for details). However, to the best of our knowledge there is no known general technique to study the Cohen-Macaulayness of fiber cones of modules that are not ideals, and this property is understood only for a few classes of modules (see for instance \cite{{Miranda},{LinP},{CLS}}). 
 
 We address this issue using generic Bourbaki ideals, and show that they allow to reduce the study of the Cohen-Macaulay property of fiber cones of modules to the case of ideals, at least in the case when enough information on the Rees algebra $\mathcal{R}(E)$ of $E$ is available. More precisely, our main result, \cref{3.5FiberCone}, is a refined version of the following theorem. 
  
  \begin{thm} \label{introFiber} \hypertarget{introFiber}{}
   Let $R$ be a Noetherian local ring, $E$ a finite $R$-module with a rank, and let $I$ be a generic Bourbaki ideal of $E$. 
    \begin{itemize}
      \item[$($a$)$] If $\mathcal{F}(E)$ is Cohen-Macaulay, then $\mathcal{F}(I)$ is Cohen-Macaulay.
      \item[$($b$)$] Assume that after a generic extension the Rees algebra $\,\mathcal{R}(E)$ is a deformation of $\,\mathcal{R}(I)$. If $\mathcal{F}(I)$ is Cohen-Macaulay, then $\mathcal{F}(E)$ is Cohen-Macaulay.
    \end{itemize}
  \end{thm}

 The deformation condition in \cref{introFiber}(b) is satisfied in particular whenever the Rees algebra of $E$ is Cohen-Macaulay. In general, the Cohen-Macaulayness of $\mathcal{F}(E)$ and of $\mathcal{R}(E)$ are not related to each other. Indeed, suppose that $R$ is Cohen-Macaulay and that $E$ isomorphic to an $R$-ideal $I$ of positive grade. Then the Rees algebra $\,\mathcal{R}(E)\,$ is isomorphic to the subalgebra
     $$\,\mathcal{R}(I) = R[It]= \oplus_{j \geq 0\,} I^j t^j\,$$ 
 of the polynomial ring $R[t]$, and is known to be Cohen-Macaulay whenever the associated graded ring  $\,\mathcal{G}(I) = \displaystyle{\oplus_{j \geq 0\,} I^j /I^{j+1}\,}$ is Cohen-Macaulay and some additional numerical conditions are satisfied \cite{{HuGrI},{IkedaTrung},{JK},{SUVDegPolyRel}}. However, one can construct perfect ideals $I$ of height two over a power series ring over a field so that $\,\mathcal{F}(I)\,$ is Cohen-Macaulay while $\,\mathcal{R}(I)\,$ is not (see the introduction of \cite{CGPU}). Moreover, for some ideals $I$ defining monomial space curves one has that $\,\mathcal{R}(I)\,$ is Cohen-Macaulay while $\,\mathcal{F}(I)\,$ is not \cite{DAnna}.  
 
 Nevertheless, in some circumstances the Cohen-Macaulay property of the Rees algebra $\,\mathcal{R}(I)\,$ of an ideal $I$ implies that of the fiber cone $\,\mathcal{F}(I)$ (see for instance \cite{{CGPU},{Jonathan}}), so it makes sense to investigate similar connections in the case of Rees algebras and fiber cones of modules as well. Our \cref{GenCGPU3.4} and \cref{GenCGPU3.1} provide classes of modules whose Rees algebra and fiber cone are both Cohen-Macaulay, generalizing previous work of Corso, Ghezzi, Polini and Ulrich on the fiber cone of an ideal (see \cite[3.1 and 3.4]{CGPU}). 
 
 A class of modules with Cohen-Macaulay fiber cone but non-Cohen-Macaulay Rees algebra is instead given in \cref{GenBM}. In this case, the deformation condition in \cref{introFiber}(b) is guaranteed by \cref{NewSUV3.7}, which is a crucial technical result in this work, as it in fact extends the applicability of generic Bourbaki ideals to the study of Rees algebras which are not necessarily Cohen-Macaulay, nor even $S_2$.

 The second main goal of this work is to study the defining ideal of Rees algebras of modules. Recall that for an $R$-module $\,E=Ra_1 + \ldots + Ra_n$ the \emph{defining ideal} of $\, \mathcal{R}(E)\,$ is the kernel of the natural homogeneous epimorphism
   \begin{eqnarray*}
      \phi \colon R[T_1, \ldots, T_n] & \longrightarrow & \mathcal{R}(E) \\ 
      T_i & \mapsto & a_i \in [\mathcal{R}(E)]_1
   \end{eqnarray*}
 
 Determining the defining ideal is usually a difficult task, but it becomes treatable for Rees algebras of ideals or modules whose free resolutions have a rich structure. Using generic Bourbaki ideals, in \cite[4.11]{SUV2003} Simis, Ulrich and Vasconcelos determined the defining ideal of $\,\mathcal{R}(E)\,$ in the case when $E$ is a module of projective dimension one with linear presentation matrix over a polynomial ring $k[X_1, \ldots, X_d]$, where $k$ is a field. Their proof ultimately relies on the fact that a generic Bourbaki ideal $I$ of $E$ has Cohen-Macaulay Rees algebra $\,\mathcal{R}(I),\,$ which allows to deduce the shape of the defining ideal of $\,\mathcal{R}(E)\,$ from that of $\,\mathcal{R}(I)$. In fact, their proof only requires that after a generic extension $\,\mathcal{R}(E)\,$ is a deformation of $\,\mathcal{R}(I)$.
 
 With a similar approach, the deformation condition of \cref{NewSUV3.7} allows us to describe the defining ideal of the Rees algebra of an \emph{almost linearly presented} module $E$ of projective dimension one over $k[X_1, \ldots, X_d]$ (see \cref{GenBM}). This condition means that all entries in a presentation matrix of $E$ are linear, except possibly those in one column, which are assumed to be homogeneous of degree $m \geq 1$. 
 
 Our result generalizes work of Boswell and Mukundan \cite[5.3]{BM} on the Rees algebra of almost linearly presented perfect ideals of height two. While this manuscript was being written, in his Ph.D. thesis \cite{Weaver} Matthew Weaver extended Boswell and Mukundan's techniques to linearly presented perfect ideals of height two over a hypersurface ring $\,\displaystyle{R = k[X_1, \ldots, X_d] /(f)}\,$, and used our methods to determine the defining ideal of the Rees algebra of linearly presented modules of projective dimension one. His work suggests potential applications to the case of Rees algebras and fiber cones of modules of projective dimension one over complete intersection rings, which include the module of K\"ahler differentials of such a ring $R$. This is particularly interesting from a geometrical perspective, since its fiber cone is the homogeneous coordinate ring of the tangential variety to the algebraic variety defined by the ring $R$.

 We now briefly describe how this paper is structured.
 
 In \cref{SecPrelim} we give the necessary background on Rees algebras and fiber cones of modules and set up the notation that will be used throughout the paper. In particular, we briefly review the construction and main properties of generic Bourbaki ideals from \cite{SUV2003}, as well as Boswell and Mukundan's construction of \emph{iterated Jacobian duals} \cite{BM}, which we will need later in \cref{SecDefEqs}. 
 
 \cref{SecDeform} contains our main technical result, namely the deformation condition of \cref{NewSUV3.7}, which is going to be crucial throughout the paper and in particular in the proofs of \cref{3.5FiberCone}, \cref{FiberType} and \cref{GenBM}.  
 
 In \cref{SecFiberCone} we study the Cohen-Macaulay property of fiber cones of modules via generic Bourbaki ideal. Our main results are \cref{3.5FiberCone}, which reduces the problem to the case of fiber cones of ideals, as well as \cref{GenCGPU3.4} and \cref{GenCGPU3.1}, which produce modules with Cohen-Macaulay fiber cones.
 
\cref{SecDefEqs} is dedicated to the study of the defining ideal of Rees algebras of modules. Besides the aforementioned \cref{GenBM} on almost linearly presented modules of projective dimension one, another key result in this section is \cref{FiberType}, which characterizes the fiber type property of a module over a standard graded $k$-algebra, where $k$ is a field.


\section{Preliminaries} \label{SecPrelim} \hypertarget{SecPrelim}{}

In this section we recall the definitions and main properties of Rees algebras and fiber cones of modules, and review the construction of generic Bourbaki ideals. 

\subsection{Rees algebras and fiber cones of modules}

 Unless otherwise specified, throughout this work, $R$ will be a Noetherian local ring and all modules will be assumed to have a rank. Recall that a finite $R$-module $E$ has a \emph{rank}, $\mathrm{rank}_{\,}E =e$, if $E \otimes_R \mathrm{Quot}(R) \cong (\mathrm{Quot}(R))^e$, or, equivalently, if $E_{\mathfrak{p}} \cong R_{\mathfrak{p}} ^e$ for all $\mathfrak{p} \in \mathrm{Ass}(R)$. 
 
 This is not a restrictive assumption. In fact, if $R$ is a domain all finite $R$-modules have a rank. Moreover, every module with a finite free resolution over any Noetherian ring has a rank.  Most importantly, torsion-free modules of rank one are isomorphic to ideals of positive grade, which is crucial for our purposes.
  
 In this setting, let $R^s \stackrel{\varphi}\longrightarrow R^n \twoheadrightarrow E$ be any presentation of $\,E =Ra_1+ \ldots +Ra_n$. Then, the natural homogeneous epimorphism
   \begin{eqnarray*}
      \phi \colon R[T_1, \ldots, T_n] & \longrightarrow & \mathcal{S}(E) \\ 
      T_i & \mapsto & a_i \in E = [\mathcal{S}(E)]_1
   \end{eqnarray*}
 onto the \emph{symmetric algebra} of $E$ induces an isomorphism 
     $$\mathcal{S}(E) \cong R[T_1, \dots ,T_n]/ \mathcal{L},$$
 where the ideal $\mathcal{L}$ is generated by linear forms $\,\ell_1, \dots, \ell_s\,$ in $\,R[T_1, \dots ,T_n]\,$ so that 
     $$ [T_1, \dots ,T_n] \cdot \varphi= [\ell_1, \dots, \ell_s].$$ 
 This definition is independent of the choice of the presentation matrix $\varphi$ (see for instance \cite[Section 1.6]{BH}). 
    
 The \emph{Rees algebra} $\mathcal{R}(E)$ of $E$ is the quotient of $\mathcal{S}(E)$ modulo its $R$-torsion submodule. In particular, 
     $$ \mathcal{R}(E) \cong R[T_1, \dots ,T_n]/ \mathcal{J} $$
 for some ideal $\mathcal{J}$, called the \emph{defining ideal} of $\mathcal{R}(E)$. Notice that by construction $\,\mathcal{J} \supseteq \mathcal{L}\,$ and the module $E$ is said to be of \emph{linear type} if equality holds, since in this case $\mathcal{J}$ is generated by linear equations. 
 
 Let $k$ be the residue field of $R$. The \emph{fiber cone} (or \emph{special fiber ring}) of $E$ is defined as 
    $$ \mathcal{F}(E) \coloneq \mathcal{R}(E) \otimes_R k $$
 (see \cite[2.3]{EHU}). It can be described as 
   $$ \mathcal{F}(E) \cong k[T_1, \dots ,T_n]/ \mathcal{I}$$
 for some ideal $\,\mathcal{I}$ in $\,k[T_1, \dots ,T_n]$. The Krull dimension $\,\ell(E) \coloneq \mathrm{dim}_{\,}\mathcal{F}(E)\,$ is called the \emph{analytic spread} of $E$ (see \cite[2.3]{EHU}) and satisfies the inequality
    $$ \,e \leq \ell(E) \leq \mathrm{dim}\,R +e-1 $$
 whenever $\,\mathrm{dim}\,R>0\,$ and $\,\mathrm{rank}\,E=e$ (see \cite[2.3]{SUV2003}).

 Similarly as for powers of an ideal $I$, one defines the power $E^j$ of a module $E$ as the $j$-th graded component of the Rees algebra $\mathcal{R}(E)$. A \emph{reduction} of $E$ is a submodule $\,U \subseteq E\,$ so that $\,E^{r+1} = U E^r\,$ for some integer $\,r \geq 0$. The least such $r$ is denoted by $r_U(E)$. A reduction $U$ of $E$ is a \emph{minimal reduction} if it is minimal with respect to inclusion and the \emph{reduction number} of $E$ is 
    $$ \,  r(E) \coloneq \mathrm{min} \, \{r _U(E) \, | \, U \mathrm{\,is \; a \; minimal \; reduction \; of\, } E\}$$
 (see \cite[2.3]{EHU}). Moreover, if $k$ is infinite then any minimal reduction of $E$ is generated by $\ell(E)$ elements, and any general $\ell(E)$ elements in $E$ generate a minimal reduction $U$ of $E$ with $\,r_U(E)=r(E)$.

\subsection{Generic Bourbaki ideals}
 
 Generic Bourbaki ideals were introduced by Simis, Ulrich and Vasconcelos in \cite{SUV2003} as a tool to study the Cohen-Macaulay property of Rees algebras of modules. Most of our technical work in this paper will consist in providing modifications of the known theory of generic Bourbaki ideals in order to extend their applicability to new situations. For this reason, we recall their construction and main properties below, referring the reader to \cite{SUV2003} for the proofs. 

 \begin{notat} \label{trueNotationBourbaki} \hypertarget{trueNotationBourbaki}{}
   \em{(\cite[3.3]{SUV2003}). Let $(R,\mathfrak{m})$ be a Noetherian local ring, $E$ a finite $R$-module with $\mathrm{rank}_{\,}E=e>0$. Let $U=Ra_1 + \dots + Ra_n$ be a submodule of $E$ for some $a_i \in E$, and consider a set of indeterminates
    \begin{displaymath}
       Z= \{ Z_{ij} \, | \, 1 \leq i \leq n, 1\leq j \leq e-1 \}.
    \end{displaymath}
   Denote $R'\coloneq R[Z]\,$ and $\,E' \coloneq E \otimes_R R'.\,$ For $\,1 \leq j \leq e-1$, let $ \, x_j= \sum_{i=1}^n Z_{ij} a_i \in E'\,$ and $\, F'=\sum_{j=1}^{e-1} R' x_j. \,$
   Also, denote $\,\displaystyle R'' = R(Z)= R[Z]_{\mathfrak{m}\,R[Z]},$ $\,\displaystyle E''=E \otimes_R R''\,$ and $\, \displaystyle F''= F' \otimes_{R'} R''$.}
 \end{notat}
 
 In the setting of \cref{trueNotationBourbaki}, the existence of generic Bourbaki ideals is guaranteed by the following result, and exploits the fundamental fact that torsion-free modules of rank one are isomorphic to ideals of positive grade.
 
 \begin{tdefn} \label[theorem]{truetdefBourbaki} \hypertarget{truetdefBourbaki}{}
   $($\cite[3.2 and 3.3]{SUV2003}$)$. Let $R$ be a Noetherian local ring, and $E$ a finite $R$-module with $\mathrm{rank}_{\,}E=e>0$, $U \subseteq E$ a submodule. Also, assume that:
         \begin{itemize}
           \item[$($i$)$] $E$ is torsion-free.
           \item[$($ii$)$] $\,E_{\mathfrak{p}}$ is free for all $\,\mathfrak{p} \in \mathrm{Spec}(R)$ with $\,\mathrm{depth}_{\,}R_\mathfrak{p} \leq 1$.
           \item[$($iii$)$] $\,\mathrm{grade}(E/U) \geq 2$.
          \end{itemize}
        Then, for $R'$, $E'$ and $F'$ as in \cref{trueNotationBourbaki}, $F'$ is a free $R'$-module of rank $e-1$ and $\,E'/F'$ is isomorphic to an $R'$-ideal $J$ with $\mathrm{grade}_{\,}J >0$. Also, $E''/F''$ is isomorphic to an $R''$-ideal $I$, called a \emph{generic Bourbaki ideal} of $E$ with respect to $U$. If $U=E$, $I$ is simply called a \emph{generic Bourbaki ideal} of $E$.
 \end{tdefn}
 
  Generic Bourbaki ideals of $E$ with respect to a submodule $U$ are essentially unique. Indeed, if $K$ is another ideal constructed as in \cref{truetdefBourbaki} using variables $Y$, then the ideals generated by $I$ and $K$ in $T=R(Z,Y)$ coincide up to multiplication by a unit in $\mathrm{Quot}(T)$, and are equal whenever $I$ and $K$ have grade at least 2 (see \cite[3.4]{SUV2003}).

  Notice that assumption (iii) in \cref{truetdefBourbaki} is automatically satisfied if $U$ is a minimal reduction of $E$. Moreover, if in this case $I \cong E''/F''$ is a generic Bourbaki ideal with respect to $U$, then the ideal $\,K \cong U''/F''$ is a minimal reduction of $I$.
  Sometimes it is possible to relate the reduction number of $E$ and the reduction number of $I$, as described in part (d) of the following theorem, which summarizes the main properties of generic Bourbaki ideals.
 
 \begin{thm} \label{trueMainBourbaki} \hypertarget{trueMainBourbaki}{}
   $($\cite[3.5]{SUV2003}$)$.
   In the setting of \cref{trueNotationBourbaki}, let $U$ be a reduction of $E$. Let $I$ be a generic Bourbaki ideal of $E$ with respect to $U$, and let $K \cong U''/F''$. Then the following statements hold.
    \begin{itemize}
      \item[$($a$)$] $\mathcal{R}(E)$ is Cohen-Macaulay if and only if $\,\mathcal{R}(I)\,$ is Cohen-Macaulay.
      \item[$($b$)$] $E$ is of linear type and $\,\mathrm{grade} \, \mathcal{R}(E)_+ \geq e\,$ if and only if $I$ is of linear type, if and only if $J$ is of linear type.
      \item[$($c$)$] If any of condition (a) or (b) hold, then $\mathcal{R}(E'')/(F) \cong \mathcal{R}(I)$ and $\,x_1, \ldots, x_{e-1}\,$ of $F$ form a regular sequence on $\mathcal{R}(E'')$.
      \item[$($d$)$] If $\,\mathcal{R}(E'')/(F) \cong \mathcal{R}(I), \,$ then $K$ is a reduction of $I$ with $r_K(I)=r_U(E)$. In this case, if in addition the residue field of $R$ is infinite and $U=E$, then $r(E)=r(I)$.
    \end{itemize}
 \end{thm}
 
 Condition (c) above says that $\mathcal{R}(E'')$ is a \emph{deformation} of $\mathcal{R}(I)$. This is in fact the key property that allows to transfer properties from $\mathcal{R}(E)$ to $\mathcal{R}(I)$ and backwards. The following result characterizes the deformation property along a Bourbaki exact sequence. 
 
  \begin{thm} \label{SUV3.11} \hypertarget{SUV3.11}{}
  $($\cite[3.11]{SUV2003}$)$.
   Let $R$ be a Noetherian ring, $E$ a finite $R$-module with $\mathrm{rank}_{\,}E=e>0$. Let $\,0 \to F \to E \to I \to 0\,$ be an exact sequence where $F$ is a free $R$-module with free basis  $x_1, \ldots, x_{e-1}$ and $I$ is an $R$-ideal. The following are equivalent.
   \begin{itemize}
       \item[$($a$)$] $\mathcal{R}(E)/(F)$ is $R$-torsion free.
       \item[$($b$)$] $\mathcal{R}(E)/(F) \cong \mathcal{R}(I)$.
       \item[$($c$)$] $\mathcal{R}(E)/(F) \cong \mathcal{R}(I)$ and $x_1, \ldots, x_{e-1}$ of $F$ form a regular sequence on $\mathcal{R}(E)$.
    \end{itemize}
  Moreover, if $I$ is of linear type, then so is $E$ and the equivalent conditions above hold.  
\end{thm}

 For our purposes, it will often be convenient to think of the rings $R'$ and $R''$ as the result of an iterative process, where at each step only $n$ variables are adjoined. This is formalized in the following notation.

 \begin{notat} \label{IterativeNotation} \hypertarget{IterativeNotation}{}
   \em{Let $R$ be a Noetherian ring, $E$ a finite $R$-module with positive rank, $U=Ra_1 + \dots + Ra_n$ a submodule of $E$ for some $a_i \in E$. Let $Z_1, \ldots Z_n$ be indeterminates, $\,\widetilde{R} \coloneq R[Z_1, \ldots, Z_n]$, $\,\widetilde{E} \coloneq E \otimes_R \widetilde{R}$, $\, \widetilde{U} \coloneq U \otimes_R \widetilde{R}$, and $\,x \coloneq \sum_{i=1}^n Z_{i} a_i \in \widetilde{U}$. If $R$ is local with maximal ideal $\mathfrak{m}$, let $\,\displaystyle{S \coloneq R(Z_1, \ldots, Z_n)= \widetilde{R}_{\mathfrak{m}\widetilde{R}}}$.} 
 \end{notat}
 
 In fact, the rings $R'$ and $R''$ as in \cref{trueNotationBourbaki} are respectively obtained from $R$ by iterating the construction of the rings $\widetilde{R}$ and $S$ as in \cref{IterativeNotation} $e-1$ times. Moreover, in \cref{trueMainBourbaki} the Cohen-Macaulay property is transferred from $\mathcal{R}(E)$ to $\mathcal{R}(I)$ and backwards using the following two results iteratively.

 \begin{thm} \label{forward} \hypertarget{forward}{}
 $($\cite[3.6 and 3.8]{SUV2003}$)$
   In the setting of \cref{IterativeNotation}, assume that $\mathrm{rank}_{\,} E =e \geq 2$ and that $E/U$ is a torsion $R$-module. Let $\,\overline{E} \coloneq \widetilde{E} / \widetilde{R}x$ and $\,\overline{\mathcal{R}} \coloneq \mathcal{R}(\widetilde{E}\,)/(x)$. Then, 
     \begin{itemize}
       \item[$($a$)$] $x$ is regular on $\mathcal{R}(\widetilde{E}\,)$.
       \item[$($b$)$] The kernel of the natural epimorphism  $\, \displaystyle \pi \colon \, \overline{\mathcal{R}} \twoheadrightarrow \mathcal{R}(\overline{E})\,$ is $K=H^0_{U\overline{\mathcal{R}}}(\overline{\mathcal{R}})$ and coincides with the $\widetilde{R}$-torsion submodule of $\overline{\mathcal{R}}$.
       \item[$($c$)$] If $U$ is a reduction of $E$ and $\mathrm{grade}_{\,}\mathcal{R}(E)_{+} \geq 2$, then $\pi$ is an isomorphism.
     \end{itemize}
 \end{thm}
 
 \begin{thm} \label{backward} \hypertarget{backward}{}
  $($\cite[3.7]{SUV2003}$)$
    In the setting of \cref{IterativeNotation}, assume that $R$ is local, that $\mathrm{rank}_{\,} E =e \geq 2$ and that $U$ is a reduction of $E$. Let $\,\overline{E}\,$ denote $\,(E \otimes_R S)/xS\,$ and $\,\overline{\mathcal{R}} \coloneq \mathcal{R}(E \otimes_R S)/(x)$. If $\,\mathcal{R}(\overline{E})\,$ satisfies $S_2$, then the natural epimorphism $\, \displaystyle \pi \colon \,\overline{\mathcal{R}} \twoheadrightarrow \mathcal{R}(\overline{E})\,$ is an isomorphism, and $x$ is regular on $\mathcal{R}(E \otimes_R S)$. In particular, $\,\mathcal{R}(E)\,$ satisfies $S_2$.
  \end{thm}

 Notice that formation of Rees algebras of finite modules commutes with flat extensions (see \cite[1.3]{EHU}). Hence, one has that $\,\mathcal{R}(\widetilde{E}) \cong \mathcal{R}(E)\otimes_R \widetilde{R},\,$ as well as $\,\mathcal{R}(E \otimes_R S) \cong \mathcal{R}(E) \otimes_R S$. Therefore, tensoring with the residue field $k$ yields isomorphisms $\,\mathcal{R}(\widetilde{E}) \otimes_R k \cong \mathcal{F}(E)\otimes_R \widetilde{R},\,$ and $\,\mathcal{F}(E \otimes_R S) \cong \mathcal{F}(E) \otimes_R S$.

 \subsection{Iterated Jacobian duals} \label{SecJac} \hypertarget{SecJac}{}
 
  When studying the defining ideal of Rees algebras, the most challenging aspect usually consists in identifying its non-linear part. In many cases of interest \cite{{Vasconcelos},{SUVjacduals},{MU},{UVeqLinPres},{Morey},{Johnson},{PU99},{BM},{KPU}}, this can be done by examining some auxiliary matrices associated with $\varphi$, namely the \emph{Jacobian dual} or the \emph{iterated Jacobian duals} of $\varphi$, introduced by Vasconcelos \cite{Vasconcelos} and by Boswell and Mukundan \cite{BM} respectively. We briefly recall these notions here, as we will use them intensively in \cref{SecDefEqs}. 
  
  Although both definitions make sense over any Noetherian ring, for our purposes we assume that $\,R=k[Y_1, \ldots, Y_d]\,$ is a standard graded polynomial ring over a field $k$. Let $S=R[T_1, \ldots, T_n]$ be bigraded, and set $\, \underline{Y}= Y_1, \ldots, Y_d$, $\,\underline{T}=T_1, \ldots, T_n$. Let
      $$\,R^{s} \stackrel{\varphi}{\longrightarrow} R^n $$
 be an $\,n \times s\,$ matrix whose entries are homogeneous of constant $\underline{Y}$-degrees $\,\delta_1, \ldots, \delta_s$ along each column and assume that $I_1(\varphi) \subseteq (\underline{Y})$.
 
 \begin{tdefn} \label[definition]{JacDual} \hypertarget{JacDual}{}
 \em{\cite{Vasconcelos}$\,$ With $R$, $S$ and $\varphi$ as above, let $M=\mathrm{coker}(\varphi)$ and let $\,\ell_1, \ldots, \ell_s\,$ be linear forms in the $T_i$ variables, generating the defining ideal of the symmetric algebra $\,\mathcal{S}(M)$. Then, 
 \begin{itemize}
     \item[(a)] There exists a $d \times s$ matrix $B(\varphi)$ whose entries are linear in the $T_i$ variables and homogeneous of constant $\underline{Y}$-degrees $\,\delta_1 -1, \ldots, \delta_s -1\,$ along each column, satisfying 
        $$ [\ell_1, \ldots, \ell_s]=[\underline{T}] \cdot \varphi = [\underline{Y}] \cdot B(\varphi).$$
     $B(\varphi)$ is called a \emph{Jacobian dual} of $\varphi$.
     \item[(b)] $B(\varphi)$ is not necessarily unique, but it is if the entries of $\varphi$ are all linear. Moreover, by Cramer's rule it follows that $\,\mathcal{L} + I_d(B(\varphi)) \subseteq \mathcal{J}$.
 \end{itemize}}   
 \end{tdefn}

For a matrix $A$, let $\,(\underline{Y} \cdot A)$ denote the ideal generated by the entries of the row vector $[\underline{Y}] \cdot A$.

 \begin{tdefn} \label{IterJacDuals} \hypertarget{IterJacDuals}{} \hypertarget{IterJacDuals}{}
  \em{(\cite[4.1, 4.2 and 4.5]{BM})$\,$
    With $R$, $S$ and $\varphi$ as above, let $B_1(\varphi)=B(\varphi)$ for some Jacobian dual $B(\varphi)$ of $\varphi$. Assume that matrices $B_j(\varphi)$ with $d$ rows have been inductively constructed for $1 \leq j \leq i$, such that each $B_j(\varphi)$ has homogeneous entries of constant $\underline{Y}$-degrees and $\underline{T}$-degrees along each column. There exists a matrix $C_i$ whose entries in $S$ are homogeneous of constant $\underline{Y}$-degrees and $\underline{T}$-degrees in each column, such that $\,B_{i+1}(\varphi) \coloneq [B_i(\varphi) \,|\, C_i]\,$ satisfies
     \begin{displaymath} \label{eqIterJac} \hypertarget{eqIterJac}{}
         (\underline{Y} \cdot B_i(\varphi)) + (I_d(B_i(\varphi)) \cap (\underline{Y})) = (\underline{Y} \cdot B_i(\varphi)) + (\underline{Y} \cdot C_i).
     \end{displaymath}
    A matrix $\,B_i(\varphi)\,$ as above is called an \emph{$i$-th iterated Jacobian dual} of $\varphi$. Moreover, for all $i \geq 1$:
   \begin{itemize}
      \item[$($a$)$] The ideal $(\underline{Y} \cdot B(\varphi))+ I_d(B_i(\varphi))$ only depends on $\varphi$.
      \item[$($b$)$]  $\, (\underline{Y} \cdot B_i(\varphi))+ I_d(B_i(\varphi))= (\underline{Y} \cdot B(\varphi))+ I_d(B_i(\varphi)) \subseteq (\underline{Y} \cdot B(\varphi)) + I_d(B_{i+1}(\varphi))$. In particular, there exists an $N>0$ so that $\,(\underline{Y} \cdot B(\varphi))+ I_d(B_i(\varphi)) = (\underline{Y} \cdot B(\varphi)) + I_d(B_{i+1}(\varphi))\,$ for all $i \geq N$.
        \item[$(c)$] $(\underline{Y} \cdot B(\varphi)) + I_d(B_i(\varphi)) \subseteq ((\underline{Y} \cdot B(\varphi)) \, \colon (\underline{Y})^i)$. 
   \end{itemize}}
 \end{tdefn}

\section{A deformation condition for the Rees algebra of a module}  \label{SecDeform} \hypertarget{SecDeform}{}

 As described in \cref{SecPrelim}, transferring properties from a module $E$ to a generic Bourbaki ideal $I$ of $E$ and backwards depends on whether the Rees algebra $\mathcal{R}(E'')$ is a deformation of $\mathcal{R}(I)$. \cref{trueMainBourbaki} and \cref{SUV3.11} show that this always occurs when $I$ is of linear type, or if $\mathcal{R}(E)$ or $\mathcal{R}(I)$ are known to be Cohen-Macaulay. On the other hand, it is interesting to find alternative conditions on $E$ or $I$ that guarantee this deformation property. 
 
 Inspired by \cite[3.7]{SUV2003} (which was stated here as \cref{backward}), the following result provides a new deformation condition, which applies to ideals and modules not necessarily of linear type and whose Rees algebras are not necessarily Cohen-Macaulay (see for instance \cref{3.5FiberCone} and \cref{GenBM}). 

 \begin{thm}\label{NewSUV3.7}
    Let $R$ be a Noetherian local ring, $E$ a finite $R$-module with $\mathrm{rank}_{\,} E =e \geq 2$ and let $\,U=Ra_1 + \dots + Ra_n\,$ be a reduction of $E$. Let $S$ and $x$ be as in \cref{IterativeNotation} and denote $\,\overline{E} \coloneq (E \otimes_R S) / Sx$. 
     
    Assume that $\,\mathrm{depth}_{\,} \mathcal{R}(\overline{E}_{\mathfrak{q}}) \geq 2\,$ for all $\,\mathfrak{q} \in \mathrm{Spec}(S)\,$ such that $\,\overline{E}_{\mathfrak{q}}\,$ is not of linear type. Then, the natural epimorphism 
       $$\, \pi \colon \, \mathcal{R}(E \otimes_R S)/(x) \twoheadrightarrow \mathcal{R}(\overline{E})\,$$ 
    is an isomorphism, and $x$ is regular on $\mathcal{R}(E \otimes_R S)$.
 \end{thm}

   \emph{Proof}. We modify the proof of \cite[3.7]{SUV2003}. Since $x$ is regular on $\mathcal{R}(E \otimes_R S)$ by \cref{forward}(a),  we only need to show that $K= \mathrm{ker}(\pi)$ is zero. In fact, we only need to prove this locally at primes $\mathfrak{q} \in \mathrm{Spec}(S)$ such that $\overline{E}_{\mathfrak{q}}$ is not of linear type. Indeed, if $\overline{E}_{\mathfrak{q}}$ is of linear type, then $\,\mathcal{R}(\overline{E}_{\mathfrak{q}}) \cong \mathcal{S}(\overline{E}_{\mathfrak{q}})\,$ is isomorphic to $\,\mathcal{S}((E \otimes_R S)_{\mathfrak{q}}) / (x)\,$ by construction, whence $K_{\mathfrak{q}}=0$. 
   

   
   Let $\overline{\mathcal{R}}$ denote $\mathcal{R}(E \otimes_R S)/(x)$ and let $\,M=(\mathfrak{m}, \mathcal{R}(E \otimes_R S)_{+})\,$ be the unique homogeneous maximal ideal of $\mathcal{R}(E \otimes_R S)$. Notice that $K \subseteq H^0_M (\overline{\mathcal{R}})$. In fact, after localizing $S$ if needed, we may assume that $K$ vanishes locally on the punctured spectrum of $S$. Hence, $K$ is annihilated by a power of $\mathfrak{m}$. Also, by \cref{forward} it follows that $K$ is annihilated by a power of $\,U \overline{\mathcal{R}}, \,$ and hence by a power of $\,E \overline{\mathcal{R}} = (\overline{\mathcal{R}})_{+}, \,$ since $E$ is integral over $U$.
 
   Thus, for all $\mathfrak{q} \in \mathrm{Spec}(S)\,$ $\,K_{\mathfrak{q}} \subseteq H^0_{M_{\mathfrak{q}}} (\overline{\mathcal{R}}_{\mathfrak{q}})$ and it suffices to show that $H^0_{M_{\mathfrak{q}}}(\overline{\mathcal{R}}_{\mathfrak{q}})=0\,$  whenever $\,\overline{E}_{\mathfrak{q}}$ is not of linear type. Consider the long exact sequence of local cohomology induced by the exact sequence 
      $$ \,0 \to K_\mathfrak{q} \to \overline{\mathcal{R}}_{\mathfrak{q}} \to \mathcal{R}(\overline{E}_{\mathfrak{q}}) \to 0\, . $$
 
   Since by assumption $\,\mathrm{depth}_{\,} \mathcal{R}(\overline{E})_{\mathfrak{q}}) \geq 2$, then $\, H^i_{M_{\mathfrak{q}}}(\overline{\mathcal{R}}_{\mathfrak{q}}) \cong H^i_{M_{\mathfrak{q}}}(K_{\mathfrak{q}}) \,$ for $i=0,1$. In particular, since $\,K_\mathfrak{q} \subseteq  H^0_{M_{\mathfrak{q}}}(\overline{\mathcal{R}}_{\mathfrak{q}}), \,$ it follows that $\,0=H^1_{M_{\mathfrak{q}}}(K_{\mathfrak{q}}) \cong H^1_{M_{\mathfrak{q}}}(\overline{\mathcal{R}}_{\mathfrak{q}})$. Therefore, the exact sequence 
      $$ 0 \to \mathcal{R}((E \otimes_R S)_\mathfrak{q})(-1) \stackrel{x}{\longrightarrow} \mathcal{R}((E \otimes_R S)_\mathfrak{q})  \longrightarrow \overline{\mathcal{R}}_\mathfrak{q} \to 0 $$
   induces the exact sequence
      $$ 0 \to H^0_{M_{\mathfrak{q}}}(\overline{\mathcal{R}}_\mathfrak{q})  \longrightarrow H^1_{M_{\mathfrak{q}}}(\mathcal{R}((E \otimes_R S)_\mathfrak{q}))(-1) \stackrel{x}{\longrightarrow} H^1_{M_{\mathfrak{q}}}(\mathcal{R}((E \otimes_R S)_{\mathfrak{q}}))  \to 0\, .$$
   Now, similarly as in \cite[3.7]{SUV2003}, one can show that $\,H^1_{M_{\mathfrak{q}}}(\mathcal{R}((E \otimes_R S)_\mathfrak{q}))$ is finitely generated, as a consequence of the graded version of the Local Duality Theorem. Therefore, by the graded version of Nakayama's Lemma, it follows that $\,H^1_{M_{\mathfrak{q}}}(\mathcal{R}((E \otimes_R S)_\mathfrak{q}))=0,\,$ whence also $\,H^0_{M_{\mathfrak{q}}}(\overline{\mathcal{R}}_\mathfrak{q})=0$.  $\,\blacksquare$

\section{Cohen-Macaulay property of fiber cones of modules} \label{SecFiberCone} \hypertarget{SecFiberCone}{}
   
 In this section we examine the Cohen-Macaulay property of fiber cones of modules. We first show that the construction of generic Bourbaki ideals allows to reduce the problem to the case of ideals, as long as the passage to a generic Bourbaki ideal induces a deformation between the Rees algebras (see \cref{3.5FiberCone}). We then provide sufficient conditions for the fiber cone of a module to be Cohen-Macaulay, generalizing known results of Corso, Ghezzi, Polini and Ulrich for the fiber cone of ideals \cite[3.1 and 3.4]{CGPU}.

 Let $I$ be a generic Bourbaki ideal of $E$. The proof of \cite[3.5]{SUV2003} (which was stated here as \cref{trueMainBourbaki}) suggests that, in order to transfer the Cohen-Macaulay property from $\, \mathcal{F}(E)\,$ to $\,\mathcal{F}(I)\,$ and backwards, the natural map 
   $$ \pi \colon \mathcal{F}(E'') \cong \mathcal{R}(E'')\otimes_R k \twoheadrightarrow \mathcal{R}(I)\otimes_R k \cong \mathcal{F}(I) $$
 needs to be an isomorphism. Hence, it suffices to provide conditions on the module $E$ or on the ideal $I$ so that this isomorphism is guaranteed. Our first goal in this direction is to prove that an analogous statement as that of \cref{forward} holds for fiber cones. This is done through the next two propositions.
 
 \begin{prop} \label{3.6FiberCone} \hypertarget{3.6FiberCone}{}
   Let $(R, \mathfrak{m}, k)$ be a Noetherian local ring, $E$ a finite $R$-module with $\mathrm{rank}_{\,} E =e \geq 2$, and let $U$ be a submodule of $E$ such that $E/U$ is torsion. In the setting of \cref{IterativeNotation}, let $L$ be the kernel of the natural epimorphism 
     $$ \pi \colon \, (\mathcal{F}(E) \otimes_R \widetilde{R})/(x) \twoheadrightarrow \mathcal{R}(\widetilde{E} / \widetilde{R}x) \otimes_R k.$$
   Then,
   \begin{itemize}
      \item[$($a$)$] $L\subseteq H^0_U((\mathcal{F}(E) \otimes_R \widetilde{R})/(x))$.
      \item[$($b$)$] If in addition $U$ is a reduction of $E$ and $\mathrm{depth}_{\,}\mathcal{F}(E)>0$, then $x$ is regular on $\,\mathcal{F}(E) \otimes_R \widetilde{R}$.
   \end{itemize}
 \end{prop}

   \emph{Proof}. Let $\,\overline{\mathcal{R}}$ denote $\,\mathcal{R}(\widetilde{E}) / (x)$. By \cref{forward}(b), there is an exact sequence 
     $$0 \to K \stackrel{\iota}{\longrightarrow} \overline{\mathcal{R}} \stackrel{\pi}\longrightarrow \mathcal{R}(\widetilde{E} / \widetilde{R}x) \to 0 $$
   where $K= H^0_{U\overline{\mathcal{R}}}(\overline{\mathcal{R}})$. Tensoring with the residue field $k$, it then follows that 
     $$ L= (\iota \otimes k) (H^0_{U\overline{\mathcal{R}}}(\overline{\mathcal{R}}) \otimes_R k) \subseteq H^0_U((\mathcal{F}(E) \otimes_R \widetilde{R})/(x)) .$$
   This proves (a). Part (b) follows from \cite{Hochster}, after noticing that $\mathrm{depth}_{\,}\mathcal{F}(E) = \mathrm{grade}_{\,}U\mathcal{F}(E)$. $\blacksquare$ \\

 \begin{prop} \label{3.8FiberCone} \hypertarget{3.8FiberCone}{}
   Let $(R, \mathfrak{m}, k)$ be a Noetherian local ring, $E$ a finite $R$-module with $\mathrm{rank}_{\,} E =e \geq 2$, and let $U$ be a reduction of $E$. With the notation of \cref{3.6FiberCone}, assume that $\mathrm{depth}_{\,}\mathcal{F}(E) \geq 2$. Then, 
    \begin{displaymath}
      \pi \colon \, (\mathcal{F}(E) \otimes_R \widetilde{R})/(x) \twoheadrightarrow \mathcal{R}(\widetilde{E} / \widetilde{R}x) \otimes_R k.
    \end{displaymath}
   is an isomorphism, and $x$ is regular on $\,\mathcal{F}(E) \otimes_R\widetilde{R}$.
 \end{prop}

   \emph{Proof}. Let $\overline{\mathcal{R}}$ denote $\,\mathcal{R}(E\otimes_R \widetilde{R})/(x)$. By \cref{3.6FiberCone} it follows that $\,L= \mathrm{ker}(\pi) = (\iota \otimes k) (H^0_{U\overline{\mathcal{R}}}(\overline{\mathcal{R}}) \otimes_R k) \subseteq H^0_U((\mathcal{F}(E) \otimes_R \widetilde{R})/(x))\,$ and that $x$ is regular on $\,\mathcal{F}(E) \otimes_R\widetilde{R}$.  
   In particular, there is an exact sequence
     $$ 0 \to (\mathcal{F}(E) \otimes_R \widetilde{R})(-1)  \stackrel{x}\longrightarrow \mathcal{F}(E) \otimes_R \widetilde{R} \longrightarrow (\mathcal{F}(E) \otimes_R \widetilde{R})/(x) \to 0.$$
   Now, notice that $\,H^1_U(\mathcal{F}(E) \otimes_R \widetilde{R})=0\,$ since 
     $$\mathrm{grade}_{\,}U \mathcal{F}(E) \otimes_R \widetilde{R}=\mathrm{grade}_{\,}E \mathcal{F}(E) \otimes_R \widetilde{R} \geq \mathrm{depth}_{\,}\mathcal{F}(E) \geq 2. $$
   Hence, the long exact sequence of local cohomology implies that $$\,H^0_U((\mathcal{F}(E) \otimes_R \widetilde{R})/(x))=0.$$ 
   Thus, $L=0$ and $\pi$ is an isomorphism. $\blacksquare$ \\
   
 By applying \cref{3.8FiberCone} repeatedly, we obtain the following useful corollary.
 \begin{cor} \label{Fiberforward} \hypertarget{Fiberforward}{}
   Let $R$ be a Noetherian local ring, and let $E$ be a finite $R$-module with $\mathrm{rank}_{\,} E =e$. In the setting of \cref{trueNotationBourbaki}, let $I$ be a generic Bourbaki ideal of $E$ with respect to a reduction $U$ of $E$. Assume that $\,\mathrm{depth}_{\,}\mathcal{F}(E) \geq e$, then the natural epimorphism 
        $$ \pi \colon \,\mathcal{F}(E'')/(F'') \twoheadrightarrow \mathcal{F}(I)\,$$
      is an isomorphism and $\,F''\mathcal{F}(E'')\,$ is generated by a regular sequence of linear forms.
 \end{cor}

 We now proceed to set up the technical framework in order for the Cohen-Macaulay property to be transferred from $\,\mathcal{F}(I)\,$ back to $\,\mathcal{F}(E)$. The key result is \cref{3.7FiberCone}, whose proof relies on the next two lemmas. \cref{filter-regular} states that $x$ is a \emph{filter-regular element} on $\,\mathcal{F}(E) \otimes_R \widetilde{R}\,$ with respect to the ideal $\,E(\mathcal{F}(E) \otimes_R \widetilde{R})$ (see for instance \cite[p. 13]{RV}).

 \begin{lemma} \label{filter-regular} \hypertarget{filter-regular}{}
   Let $R$ be a Noetherian local ring, $E$ a finite $R$-module with $\mathrm{rank}_{\,} E =e \geq 2$, and let $U$ be a reduction of $E$. Then, in the setting of \cref{IterativeNotation}, 
      $$ \mathrm{Supp}_{\mathcal{F}(E) \otimes_R \widetilde{R}\,}(0 \colon_{\!\mathcal{F}(E) \otimes_R \widetilde{R}\,} x) \subseteq V(E(\mathcal{F}(E) \otimes_R \widetilde{R})).$$
 \end{lemma}

   \emph{Proof}. Let $\, \mathfrak{q} \in \mathrm{Spec}(\mathcal{F}(E) \otimes_R \widetilde{R}) \setminus V(E(\mathcal{F}(E) \otimes_R \widetilde{R}))\,$ and let $\,\mathfrak{p}=\mathfrak{q} \cap \mathcal{F}(E)$. Since $\,U=R_{\,} a_1 + \ldots +R_{\,} a_n$ is a reduction of $E$, it follows that $\,\mathfrak{q} \nsupseteq U(\mathcal{F}(E) \otimes_R \widetilde{R})$. Hence, $\,\mathfrak{p} \nsupseteq U\mathcal{F}(E),\,$ which means that at least one of the $a_i$ is a unit in $\,\mathcal{F}(E)_{\mathfrak{p}}$. Therefore, $\,x=\sum_{i=1}^n Z_{i} a_i\,$ is a nonzerodivisor in $\,\mathcal{F}(E)_{\mathfrak{p}}[Z_1, \ldots, Z_n]$, hence also in its further localization $\,(\mathcal{F}(E)[Z_1, \ldots, Z_n])_{\mathfrak{q}} = (\mathcal{F}(E) \otimes_R \widetilde{R})_{\mathfrak{q}}$. Thus, $\,(0 \colon_{\!\mathcal{F}(E) \otimes_R \widetilde{R}\,} x)_{\mathfrak{q}} =0$. $\blacksquare$\\

 \begin{lemma} \label{homework} \hypertarget{homework}{}
   Let $R$ be a positively graded Noetherian ring with $R_0$ local, and let $x$ be a homogeneous non-unit element of $R$. Let $M$ be a finite graded $R$-module, and assume that $\mathrm{dim}_{\,}(0 \colon_{\!\!M \,}x) < \mathrm{depth}_{\,}(M/xM)$. Then, $x$ is a nonzerodivisor on $M$.
 \end{lemma}
 
   \emph{Proof}. It suffices to show that $\,H^0_{(x)}(M)=0$,
   which would follow from Nakayama's Lemma once we prove that $\,H^0_{(x)}(M) / x H^0_{(x)}(M)=0$. For this, it suffices to show that $\, \mathrm{Ass}(H^0_{(x)}(M) / x H^0_{(x)}(M)) = \emptyset$. 
   
   Consider the short exact sequences 
    \begin{equation} \label{sescolon}
      0 \to 0 \colon_{\!\!M\,} x \to M \to M/(0 \colon_{\!\!M\,} x) \to 0
    \end{equation}
   and 
    \begin{equation} \label{sesxM}
      0 \to M/(0 \colon_{\!\!M\,} x) \stackrel{x}{\longrightarrow} M \to M/xM \to 0.
    \end{equation}
   Since $\,(0 \colon_{\!\!M\,} x) = H^0_{(x)}(0 \colon_{\!\!M\,} x),\,$ it follows that $\,H^1_{(x)}(0 \colon_{\!\!M\,} x)=0$. Hence, the long exact sequence of local cohomology induced by~(\ref{sescolon}) implies that $\, H^0_{(x)}(M)\,$ surjects onto $\,H^0_{(x)}(M/(0 \colon_{\!\!M\,} x))$. Therefore, the long exact sequence of local cohomology induced by~(\ref{sesxM}) 
     $$ 0 \to H^0_{(x)}(M/(0 \colon_{\!\!M\,} x)) \stackrel{x}{\longrightarrow} H^0_{(x)}(M) \to H^0_{(x)}(M/xM) $$
   in turn induces an exact sequence
     $$H^0_{(x)}(M) \stackrel{x}{\longrightarrow} H^0_{(x)}(M) \to H^0_{(x)}(M/xM) \subseteq M/xM.$$
   In particular,  $\, H^0_{(x)}(M) / x H^0_{(x)}(M) \,$ embeds into $\,M/xM$.
   
   Also, notice that $(0 \colon_{\!} x)=0$ if and only if $\,H^0_{(x)}(M)=0$, hence if and only if $\,H^0_{(x)}(M) / x H^0_{(x)}(M)=0\,$ by Nakayama's Lemma. Therefore, $\, \mathrm{Supp}(0 \colon_{\!\!M\,} x)$ $= \mathrm{Supp}(H^0_{(x)}(M) / x H^0_{(x)}(M))$. Hence, if there exists some $\, \mathfrak{p} \in \mathrm{Ass}(H^0_{(x)}(M) / x H^0_{(x)}(M))$, then $\,\mathrm{dim}(R/\mathfrak{p}) \leq \mathrm{dim}(0 \colon_{\!\!M\,} x)$. On the other hand, since $\, \mathrm{Ass}(H^0_{(x)}(M) / x H^0_{(x)}(M)) \subseteq \mathrm{Ass}(M/xM), \,$ we also have that $\, \mathrm{dim}(R/\mathfrak{p}) \geq \mathrm{depth}(M/xM)$. But then
    $$\,\mathrm{dim}(0 \colon_{\!\!M\,} x) = \mathrm{dim}(H^0_{(x)}(M) / x H^0_{(x)}(M)) \geq \mathrm{depth}(M/xM), \,$$ 
   which contradicts the assumption. So, it must be that $\, \mathrm{Ass}(H^0_{(x)}(M) / x H^0_{(x)}(M)) = \emptyset$, as we wanted to prove. $\blacksquare$ \\

 \begin{thm} \label{3.7FiberCone} \hypertarget{3.7FiberCone}{}
   In the setting of \cref{IterativeNotation}, assume that $R$ is local and that $\mathrm{rank}_{\,} E =e \geq 2$. Let $\,U$ be a reduction of $E$, and denote $\,\overline{E} \coloneq (E \otimes_R S)/Sx$. Assume that one of the two following conditions hold:
    \begin{itemize}
      \item[$($i$)$] $\,\mathcal{R}(\overline{E})\,$ satisfies $S_2$, or
       \item[$($ii$)$] $\mathrm{depth}_{\,} \mathcal{R}(\overline{E}_{\mathfrak{q}}) \geq 2$ for all $\mathfrak{q} \in \mathrm{Spec}(S)$ such that $\,\overline{E}_{\mathfrak{q}}$ is not of linear type. 
     \end{itemize}
   Then, the natural epimorphism 
      $$ \pi \colon \,\mathcal{F}(E \otimes_R S)/(x) \twoheadrightarrow \mathcal{F}(\overline{E})\,$$ 
   is an isomorphism. Moreover, $x$ is regular on $\,\mathcal{F}(E \otimes_R S)\,$ if $\,\mathrm{depth}_{\,} \mathcal{F}(\overline{E})>0$.
 \end{thm}

   \emph{Proof}. Assumption (i) and \cref{backward} together imply that the natural epimorphism 
     $$\pi \colon \,\mathcal{R}(E\otimes_R S)/(x) \twoheadrightarrow \mathcal{R}(\overline{E}) $$
   is an isomorphism. The same conclusion holds if assumption (ii) is satisfied, thanks to \cref{NewSUV3.7}. Hence, $\, \displaystyle \pi \colon \,\mathcal{F}(E \otimes_R S)/(x) \twoheadrightarrow \mathcal{F}(\overline{E})\,$ is an isomorphism as well. In particular, if in addition $\,\mathrm{depth}_{\,} \mathcal{F}(\overline{E})>0\,$ then $\,\mathrm{depth}_{\,}(\mathcal{F}(E \otimes_R S)/(x)) >0$. Moreover, by \cref{filter-regular} we know that $\,(0 \colon_{\mathcal{F}(E\otimes_R S)}\, x)$ is an Artinian $\mathcal{F}(E\otimes_R S)$-module. Hence, $x$ is regular thanks to \cref{homework}. $\blacksquare$ \\

  By applying \cref{3.7FiberCone} repeatedly, we obtain the following corollary.
 \begin{cor} \label{Fiberbackward} \hypertarget{Fiberbackward}{}
   Let $R$ be a Noetherian local ring, and let $E$ be a finite $R$-module with $\mathrm{rank}_{\,} E =e$. In the setting of \cref{trueNotationBourbaki} $I$ be a generic Bourbaki ideal of $E$ with respect to a reduction $U$ of $E$. 
    \begin{itemize}
      \item[$($a$)$] Assume that either $\, \mathcal{R}(I)$ is $S_2$, or $\,\mathrm{depth}_{\,} \mathcal{R}(I_{\mathfrak{q}}) \geq 2\,$ for all $\mathfrak{q} \in \mathrm{Spec}(R'')$ so that $I_{\mathfrak{q}}$ is not of linear type. Then, the natural epimorphism 
        $$ \pi \colon \,\mathcal{F}(E'')/(F'') \twoheadrightarrow \mathcal{F}(I)\,$$
      is an isomorphism.
      \item[$($b$)$] If in addition $\mathcal{F}(I)$ is Cohen-Macaulay, then $\,F''\mathcal{F}(E'')\,$ is generated by a regular sequence of linear forms. 
    \end{itemize}
 \end{cor}

 \emph{Proof}. Notice that the assumptions in (a) imply that the assumptions (i) or (ii) in \cref{3.7FiberCone} are satisfied at each iteration, thanks to \cref{backward} or \cref{NewSUV3.7} respectively. Hence, $\,\mathcal{R}(E'')/(F'') \cong \mathcal{R}(I),\,$ and $F''\mathcal{R}(E'')$ is generated by a regular sequence on $\mathcal{R}(E'')$. Hence, by iteration of \cref{3.7FiberCone}, it follows that $\,\mathcal{F}(E'')/(F'') \cong \mathcal{F}(I)$. Now, if furthermore $\,\mathcal{F}(I)$ is Cohen-Macaulay, then the proof of \cref{3.7FiberCone} implies that also $\,F''\mathcal{F}(E'')$ is generated by a regular sequence on $\mathcal{F}(E'')$. That the generators of $\,F''\mathcal{F}(E'')$ are linear forms in $\mathcal{F}(E'')$ is clear by construction. $\blacksquare$ \\

   

 We are now ready to state and prove our main result. 

 \begin{thm} \label{3.5FiberCone} \hypertarget{3.5FiberCone}{}
    Let $R$ be a Noetherian local ring, $E$ a finite $R$-module with $\mathrm{rank}_{\,} E =e$, $U$ a reduction of $E$. Let $I$ be a generic Bourbaki ideal of $E$ with respect to $U$. 
    \begin{itemize}
      \item[$($a$)$] If $\mathcal{F}(E)$ is Cohen-Macaulay, then $\mathcal{F}(I)$ is Cohen-Macaulay.
      \item[$($b$)$] Assume that either $\, \mathcal{R}(I)$ is $S_2$, or $\,\mathrm{depth}_{\,} \mathcal{R}(I_{\mathfrak{q}}) \geq 2\,$ for all $\mathfrak{q} \in \mathrm{Spec}(R'')$ so that $I_{\mathfrak{q}}$ is not of linear type. If $\mathcal{F}(I)$ is Cohen-Macaulay, then $\mathcal{F}(E)$ is Cohen-Macaulay.
    \end{itemize}
 \end{thm}
 
    \emph{Proof}. We may assume that $e\geq 2$. Suppose that $\,\mathcal{F}(E)\,$ is Cohen-Macaulay, then $\,\mathrm{depth}_{\,} \mathcal{F}(E) =\ell(E) \geq e$. Hence, by \cref{Fiberforward} it follows that $\,\mathrm{depth}_{\,} \mathcal{F}(I) =\ell(E) -e+1$. The latter equals $\ell(I)$ by \cite[3.10]{SUV2003}, hence $\,\mathcal{F}(I)\,$ is Cohen-Macaulay. 

    Conversely, if the assumptions in (b) hold, by \cref{Fiberforward} it follows that $\,\mathrm{depth}_{\,} \mathcal{F}(E) = \ell(I)+e-1 =\ell(E)$. Therefore, $\,\mathcal{F}(E)\,$ is Cohen-Macaulay. $\blacksquare$ \\

 From the proofs of \cref{Fiberforward} and \cref{3.5FiberCone} it follows that $\,\mathcal{F}(E)\,$ is Cohen-Macaulay whenever $\,\mathcal{F}(I)\,$ is Cohen-Macaulay and $\,\mathcal{R}(E'')\,$ is a deformation of $\,\mathcal{R}(I)$. In particular, finding conditions other than those in \cref{backward} or \cref{NewSUV3.7} to guarantee this deformation property for the Rees algebras would provide alternative versions of \cref{3.5FiberCone}. In fact, in order to transfer the Cohen-Macaulay property from $\,\mathcal{F}(E)\,$ to $\,\mathcal{F}(I)\,$ and backwards one would only need $\,\mathcal{F}(E'')\,$ to be a deformation of $\,\mathcal{F}(I)$. Finding conditions for this to occur without any prior knowledge of the Rees algebra would potentially allow to use generic Bourbaki ideals also in the case when the Cohen-Macaulayness of $\,\mathcal{F}(I)\,$ is possibly unrelated to that of $\,\mathcal{R}(I)\,$ (see for instance \cite{{Shah},{D'CruzRV},{CGPU},{JV1},{JV2},{CPV},{Viet08},{Jonathan}}).

 \subsection{Modules with Cohen-Macaulay fiber cone} \label{SecCMFiberE} \hypertarget{SecCMFiberE}{}
 
 \cref{3.5FiberCone} above implies that the fiber cone $\,\mathcal{F}(E)\,$ of $E$ is Cohen-Macaulay whenever both the Rees algebras $\,\mathcal{R}(I)\,$ and the fiber cone $\,\mathcal{F}(I)\,$ of a generic Bourbaki ideal $I$ of $E$ are Cohen-Macaulay. The goal of this section is to provide specific classes of modules with this property. A class of modules with Cohen-Macaulay fiber cone and non-Cohen-Macaulay Rees algebra will be provided later in \cref{GenBM}.

 Our first result regards modules of projective dimension one, and extends a known result proved by Corso, Ghezzi, Polini and Ulrich for perfect ideals of height two (see \cite[3.4]{CGPU}). Recall that a module $E$ of rank $e$ satisfies \emph{condition $G_s$} if $\,\mu(E_{\mathfrak{p}}) \leq \mathrm{dim}R_\mathfrak{p} -e +1\,$ for every $\mathfrak{p} \in \mathrm{Spec}(R)$ with $1\leq \mathrm{dim}_{\,}R_{\mathfrak{p}} \leq s-1$. Moreover, by \cite[3.2]{SUV2003} if $E$ satisfies $G_s$ then so does a generic Bourbaki ideal $I$ of $E$.
 
 \begin{thm} \label{GenCGPU3.4} \hypertarget{GenCGPU3.4}{}
    Let $R$ be a local Cohen-Macaulay ring, and let $E$ a finite, torsion-free $R$-module with $\, \mathrm{projdim}(E)=1$, with $\ell(E)=\ell$. Assume that $E$ satisfies $G_{\ell-e+1}$. If $\,\mathcal{R}(E)$ is Cohen-Macaulay, then $\,\mathcal{F}(E)$ is Cohen-Macaulay.
 \end{thm}
 
    \emph{Proof}. Since $E$ is a torsion-free module of projective dimension one which satisfies $G_{\ell-e+1}$, then $E$ admits a generic Bourbaki ideal $I$, which is perfect of height 2 (see for instance the proof of \cite[4.7]{SUV2003}). If $e=1$, then the conclusion follows from \cite[3.4]{CGPU}. Otherwise, notice that $\, \mathcal{R}(I)\,$ is Cohen-Macaulay by \cref{trueMainBourbaki}, whence $\, \mathcal{F}(I)\,$ is Cohen-Macaulay by \cite[3.4]{CGPU}. Hence, $\, \mathcal{F}(E)$ is Cohen-Macaulay by \cref{3.5FiberCone}. $\blacksquare$

 \begin{cor}
    Let $R$ be a local Cohen-Macaulay ring, and let $E$ a finite, torsion-free $R$-module with $\, \mathrm{projdim}\,E=1$, with $\ell(E)=\ell$. Let $n=\mu(E)$ and let 
    \begin{displaymath}
        0 \to R^{\,n-e} \stackrel{\varphi}{\longrightarrow} R^n \to E \to 0 
     \end{displaymath}
    be a minimal free resolution of $E$. Assume that $E$ satisfies $G_{\ell-e+1}$ and that one of the following equivalent conditions hold.
    \begin{itemize}
      \item[$($i$)$] $\, r(E) \leq \ell-e$.
      \item[$($ii$)$] $\, r(E_{\mathfrak{p}}) \leq \ell-e\,$ for every prime $\,\mathfrak{p}$ with $\, \mathrm{dim}_{\,}R_{\mathfrak{p}}= \ell(E_{\mathfrak{p}})-e+1 = \ell-e+1$.
      \item[$($iii$)$] After elementary row operations, $\,I_{n-\ell}(\varphi)$ is generated by the maximal minors of the last $n-\ell$ rows of $\varphi$.
    \end{itemize}
    Then, $\, \mathcal{F}(E)$ is Cohen-Macaulay. 
 \end{cor}
 
    \emph{Proof}. By \cite[4.7]{SUV2003}, each of the conditions (i)-(iii) is equivalent to $\, \mathcal{R}(E)$ being Cohen-Macaulay. Hence, the conclusion follows from \cref{GenCGPU3.4}. $\blacksquare$ \\

 The following result was proved for fiber cones of ideals in \cite[3.1]{CGPU}.
     
 \begin{thm} \label{GenCGPU3.1} \hypertarget{GenCGPU3.1}{}
    Let $R$ be a Cohen-Macaulay local ring with infinite residue field. Let $E$ be a finite, torsion-free $R$-module with $\,\mathrm{rank}_{\,}E =e>0,\,$ $\,\ell(E)=\ell\,$ and $r(E)=r$. Assume that $E$ satisfies $G_{\ell-e+1}$ and $\ell-e+1 \geq 2$. Let $I$ be a generic Bourbaki ideal of $E$, and $\,g=\mathrm{ht}(I)$. Suppose that one of the following conditions holds.
     \begin{itemize}
        \item[$($i$)$] If $\mu(E) \geq \ell+2$, then $\,\mathcal{F}(E)$ has at most two homogeneous generating relations in degrees $\,\leq \mathrm{max} \{\, r, \ell-e-g+1\,\}$.
        \item[$($ii$)$] If $\mu(E) = \ell+1$, then $\,\mathcal{F}(E)$ has at most two homogeneous generating relations in degrees $\, \leq \ell-e-g+1$.
     \end{itemize}
    If $\, \mathcal{R}(E)\,$ is Cohen-Macaulay, then $\, \mathcal{F}(E)\,$ is Cohen-Macaulay.
 \end{thm}

    \emph{Proof}. By \cref{truetdefBourbaki}, $E$ admits a generic Bourbaki ideal $I$ with $\,\mathrm{ht}_{\,}I = g \geq 1, \,$ $\,\mu(I)=\mu(E)-e+1\,$ and $\,r(I) \leq r, \,$ which satisfies $G_{\,\ell-e+1}$, i.e. $G_{\,\ell(I)}$ (see \cite[3.10]{SUV2003}). If $e=1$ the conclusion follows from \cite[3.1]{CGPU}, since $\,\mathcal{G}(I)\,$ is Cohen-Macaulay whenever $\,\mathcal{R}(I)\,$ is by \cite[Proposition 1.1]{HuGrI}. So we may assume that $e \geq 2$. We show that both $\,\mathcal{R}(I)\,$ and $\,\mathcal{F}(I)\,$ are Cohen-Macaulay, whence $\,\mathcal{F}(E)\,$ is Cohen-Macaulay by \cref{3.5FiberCone}.
   
    Since by assumption $\,\mathcal{R}(E)\,$ is Cohen-Macaulay, then $\,\mathcal{R}(I)\,$ is Cohen-Macaulay by \cref{trueMainBourbaki}. Hence also the associated graded ring $\,\mathcal{G}(I)\,$ is Cohen-Macaulay by \cite[Proposition 1.1]{HuGrI}. Also, by \cref{Fiberbackward} there is a homogeneous isomorphism $\, \mathcal{F}(E'')/(F'') \cong \mathcal{F}(I)$. Therefore, if condition (i) holds, then whenever $\mu(I)\geq \ell+2-e+1 =\ell(I)+2\,$ it follows that $\, \mathcal{F}(I)$ has at most two homogeneous generating relations in degrees at most $\, \mathrm{max} \{\, r, \ell-e-g+1\,\}, \,$ hence in degrees at most $\, \mathrm{max} \{\, r(I), \ell(I)-g\,\}$. Similarly, in the situation of assumption (ii), whenever $\,\mu(I) = \ell(I)+1\,$ it follows that $\,\mathcal{F}(I)$ has at most two homogeneous generating relations in degrees at most $\, \ell-e-g+1=\ell(I)-g$. Hence, $\,\mathcal{F}(I)\,$ is Cohen-Macaulay by \cite[3.1]{CGPU}. $\blacksquare$ \\

 In particular, with some additional assumptions on the module $E$, we can give more explicit sufficient conditions for $\, \mathcal{F}(E)\,$ to be Cohen-Macaulay. The next two corollaries exploit results on the Cohen-Macaulay property of Rees rings from \cite{myCMRees} and recover \cite[2.9]{CGPU} in the case when $E$ is an ideal of grade at least two.
 
 Recall that a module $E$ is called \emph{orientable} if $E$ has a rank $\,e >0\,$ and $\, (\bigwedge^{e} E)^{\ast \ast} \cong R, \,$ where $(-)^{\ast}$ denotes the functor $\mathrm{Hom}_R (-, R))$. 

 \begin{cor} 
    Let $R$ be a local Gorenstein ring of dimension $d$ with infinite residue field. Let $E$ be a finite, torsion-free, orientable $R$-module, with $\mathrm{rank}E=e>0$ and $\ell(E)=\ell$. Let $g$ be the height of a generic Bourbaki ideal of $E$, and assume that the following conditions hold.
     \begin{itemize}
       \item[$($a$)$] $E$ satisfies $G_{\,\ell-e+1}$. 
       \item[$($b$)$] $r(E) \leq k$ for some integer $1 \leq k \leq \ell-e$. 
       \item[$($c$)$] $\displaystyle{
             \mathrm{depth}_{\,}E^j \geq \Big\lbrace \begin{array}{cc}
                       d-g-j+2\, \; \; & \mathrm{for} \; \; 1 \leq j \leq \ell-e-k-g+1 \\
                       d-\ell+e+k-j\, \; & \mathrm{for} \; \; \ell-e-k-g+2 \leq j \leq k \\
                   \end{array}}$
       \item[$($d$)$] If $g = 2$, $\,\mathrm{Ext}_{R_{\mathfrak{p}}}^{\,j+1}(E_{\mathfrak{p}}^j , R_{\mathfrak{p}}) =0\,$ for $\,\ell-e-k \leq j \leq \ell-e-3\,$ and for all $\,\mathfrak{p} \in \mathrm{Spec}(R)$ with $\, \mathrm{dim} R_{\mathfrak{p}} = \ell-e\, $ such that $E_p$ is not free.
     \end{itemize}
   Assume furthermore that one of the following two conditions holds.
    \begin{itemize}
        \item[$($i$)$] If $\mu(E) \geq \ell+2$, then $\,\mathcal{F}(E)$ has at most two homogeneous generating relations in degrees $\,\leq \mathrm{max} \{\, r, \ell-e-g+1\,\}$.
        \item[$($ii$)$] If $\mu(E) = \ell+1$, then $\,\mathcal{F}(E)$ has at most two homogeneous generating relations in degrees $\, \leq \ell-e-g+1 \,$.
    \end{itemize}
  Then, $\,\mathcal{F}(E)\,$ is Cohen-Macaulay.
 \end{cor}
 
    \emph{Proof}. Assumptions (a)-(d) together imply that $\,\mathcal{R}(E)$ is Cohen-Macaulay, thanks to \cite[4.3]{myCMRees}. Hence, if either condition (i) or (ii) hold it follows that $\,\mathcal{F}(E)\,$ is Cohen-Macaulay by \cref{GenCGPU3.1}. $\blacksquare$ \\

 Recall that an $R$-module $E$ is called an \emph{ideal module} if $E  \neq 0$ is finitely generated, torsion-free and so that $E^{**}$ is free, where ${-}^*$ denotes the functor $\mathrm{Hom}_R(-, R)$. Equivalently, $E$ is an ideal module if and only if $E$ embeds into a finite free module $G$ with $\,\mathrm{grade}(G/E) \geq 2$ (see \cite[5.1]{SUV2003}). 
 
 \begin{cor} 
   Let $R$ be a local Cohen-Macaulay ring, and let $E$ be an ideal module with $\mathrm{rank}_{\,}E=e$ and $\ell(E)=\ell$. Assume that the following conditions hold.
   \begin{itemize}
       \item[$($a$)$] $r(E) \leq k$, where $k$ is an integer such that $1 \leq k \leq \ell -e$.
       \item[$($b$)$] $E$ is free locally in codimension $\,\ell-e- \mathrm{min} \{2, k\},\,$ and satisfies $G_{\,\ell-e+1}$.
       \item[$($c$)$]$\mathrm{depth}(E^j) \geq d-\ell+e+k-j\,$ for $1 \leq j \leq k$.
   \end{itemize}
   Assume furthermore that one of the following two conditions holds.
    \begin{itemize}
        \item[$($i$)$] If $\mu(E) \geq \ell+2$, then $\,\mathcal{F}(E)$ has at most two homogeneous generating relations in degrees $\,\leq \mathrm{max} \{\, r, \ell-e-g+1\,\}$.
        \item[$($ii$)$] If $\mu(E) = \ell+1$, then $\,\mathcal{F}(E)$ has at most two homogeneous generating relations in degrees $\, \leq \ell-e-g+1 \,$.
    \end{itemize}
   Then, $\,\mathcal{F}(E)\,$ is Cohen-Macaulay.
 \end{cor}
 
   \emph{Proof}. From assumptions (a)-(c) it follows that $\,\mathcal{R}(E)$ is Cohen-Macaulay, thanks to \cite[4.10]{myCMRees}. Hence, if either condition (i) or (ii) hold, then $\,\mathcal{F}(E)\,$ is Cohen-Macaulay by \cref{GenCGPU3.1}. $\blacksquare$

 \section{Defining ideal of Rees algebras}
  \label{SecDefEqs} \hypertarget{SecDefEqs}{}

 In this section we use generic Bourbaki ideals to understand the defining ideal of Rees algebras of modules. The idea is not completely novel, as it appears in the proof of \cite[4.11]{SUV2003}, where the authors determine the defining ideal for the Rees algebra of a module $E$ of projective dimension one having a linear presentation. In their case, the Rees algebra of a generic Bourbaki ideal $I$ of $E$ is Cohen-Macaulay, however the latter condition will not be guaranteed nor required in the situations considered in this section.
 
 The first key observation is that it is always possible to relate a presentation matrix of $E$ to a presentation matrix of a generic Bourbaki ideal $I$ of $E$.

 \begin{rmk} \label{BourbakiPres} \hypertarget{BourbakiPres}{}
     \em{(See also \cite[p.617]{SUV2003}). Let $(R, \mathfrak{m}, k)$ be a Noetherian local ring and let $E$ be a finite $R$-module with a minimal presentation $\, \displaystyle R^{s} \stackrel{\varphi}{\longrightarrow} R^n \to E \to 0\,$.
     \begin{itemize}
         \item[(a)]  With $Z$ and $x_j$ as in \cref{trueNotationBourbaki}, by possibly multiplying $\varphi$ from the left by an invertible matrix with coefficients in $k(Z)$, we may assume that $\varphi\,$ presents $E''$ with respect to a minimal generating set of the form ${x_1, \ldots, x_{e-1}, a_e, \ldots, a_n}$. Then, $\, \displaystyle \varphi = \left[ \begin{array}{c}
          A \\
          \hline
          \psi \\
    \end{array} \right]$, 
    where $A$ and $\psi$ are submatrices of size $(e-1) \times s$ and $(n-e+1) \times s$ respectively. By construction, $\psi$ is a presentation of $I$, and is minimal since $\mu(I)=\mu(E)-e+1=n-e+1$.
         \item[(b)] Assume that $\,R=S_{\mathfrak{m}}$, where $S$ is a standard graded algebra over a field and $\mathfrak{m}$ is its unique homogeneous maximal ideal. If the entries of $\varphi$ are homogeneous polynomials of constant degrees $\delta_1, \ldots, \delta_s$ along each column, then the entries of $\psi$ are homogeneous polynomials of constant degrees $\delta_1, \ldots, \delta_s$ along each column.
     \end{itemize}}
 \end{rmk}

Let $R$ be a standard graded algebra over a field $k$ and let $E$ be a finite $R$-module. Then, the fiber cone $\mathcal{F}(E)$ of $E$ has a particularly useful description as a subring of a polynomial ring over $k$, which is summarized in the following remark. 
 
 \begin{rmk} \label{bigrading} \hypertarget{bigrading}{}
   \em{Let $R$ be a standard graded algebra over a field $k$ and homogeneous maximal ideal $\mathfrak{m}$. Let $\,E=Ra_1 + \ldots +Ra_n\,$ be finite $R$-module
   minimally generated by elements of the same degree. On the polynomial ring $\,S = R[T_1, \ldots T_n]$ define a bigrading by setting $\, \mathrm{deg}\,R_i =  (i,0)\,$ and $\,\mathrm{deg}\,T_i = (0,1)\,$. Then, the Rees algebra $\,\displaystyle{\mathcal{R}(E) \cong S/ \mathcal{J}} \,$ has a natural bigraded structure induced by the bigrading on $S$ and is generated in degrees $(0,1)$. Moreover, $\,\displaystyle{\mathfrak{m}\, \mathcal{R}(E) \cong [\mathcal{R}(E)]_{(>0, -)}}\,$. Hence, the fiber cone $\mathcal{F}(E)$ satisfies
      $$ \mathcal{F}(E) \,\cong \,\mathcal{R}(E) / \mathfrak{m}\mathcal{R}(E) \,\cong \,[\mathcal{R}(E)]_{(0,-)} \subseteq k[T_1, \ldots, T_n].$$
   As a consequence, the homogeneous epimorphism 
      $$ k[T_1, \ldots, T_n] = S \otimes_R k \twoheadrightarrow \mathcal{R}(E) \otimes_R k =  \mathcal{F}(E)$$
   has kernel $\, \mathcal{I} \cong \mathcal{J}_{(0,-)}$.}
 \end{rmk}
 
 Notice also that for an $R$-module $E$ as in \cref{bigrading} the defining ideal $\mathcal{L}$ of the symmetric algebra $\mathcal{S}(E)$ satisfies $\, \displaystyle{\mathcal{L} = \mathcal{J}_{(-,1)} }$. In particular, if $\mathcal{I}R[T_1, \ldots, T_n]$ denotes the extension of $\mathcal{I}$ to the ring $R[T_1, \ldots, T_n]$, then 
    $$\,\mathcal{J} \supseteq \mathcal{L} + \mathcal{I}R[T_1, \ldots, T_n].$$ 
 $E$ is said to be of \emph{fiber type} if the latter inclusion is an equality, or equivalently, if $\mathcal{J}$ is generated in bidegrees $(-,1)$ and $(0,-)$. The following theorem characterizes the fiber type property of modules.

 \begin{thm} \label{FiberType} \hypertarget{FiberType}{}
    Let $R=S_{\mathfrak{m}}$ where $S$ is a standard graded algebra over a field with unique homogeneous maximal ideal $\mathfrak{m}$. Let $E$ be a finite $R$-module with $\mathrm{rank}\,E=e\geq 2$, minimally generated by elements $\,a_1, \ldots, a_n\,$ that are images in $R$ of homogeneous elements of the same degree in $S$. Assume that $E$ is torsion-free and that $E_{\mathfrak{p}}$ is free for all $\mathfrak{p} \in \mathrm{Spec}\,R$ with $\mathrm{depth}\,R \leq 1$, and let $I$ be a generic Bourbaki ideal of $E$ constructed with respect to the generators $a_1, \ldots, a_n$.
    
    Assume that one of the following conditions hold. 
    \begin{itemize}
        \item[$($i$)$]  $\mathcal{R}(I)$ satisfies $S_2$; or
        \item[$($ii$)$] $\mathrm{depth}_{\,} \mathcal{R}(I_{\mathfrak{q}}) \geq 2\,$ for all $\mathfrak{q} \in \mathrm{Spec}(R'')$ so that $I_{\mathfrak{q}}$ is not of linear type. 
    \end{itemize}
    Then, $E$ is of fiber type if and only if $I$ is of fiber type.
 \end{thm}
 
    \emph{Proof}. Let $ \,R[T_1, \ldots, T_n] \twoheadrightarrow \mathcal{R}(E) \,$ be the natural epimorphism mapping $T_i$ to $a_i$ for all $i$. As in \cref{trueNotationBourbaki}, for $\, 1 \leq j \leq e-1\,$ let $\, \displaystyle{x_j = \sum_{i=1}^{n} Z_{ij}a_i}$. For every $j\,$ let $\, \displaystyle{X_j = \sum_{i=1}^{n} Z_{ij}T_i \subseteq R''[T_1, \ldots, T_n] \,}$ and notice that $X_j$ is mapped to $x_j$ via the natural epimorphism $\, \displaystyle{R''[T_1, \ldots, T_n] \twoheadrightarrow \mathcal{R}(E'')}$. Let $\varphi$ be a minimal presentation of $E$ with respect to the generators $\,a_1, \ldots, a_n,\,$ let $\psi$ be a minimal presentation of $I$ constructed as in \cref{BourbakiPres}, and use these presentations to construct the symmetric algebras $\mathcal{S}(E)$ and $\mathcal{S}(I)$ respectively. Let $\mathcal{L}_E$, $\mathcal{J}_E$ and $\mathcal{I}_E$ denote the defining ideals of $\mathcal{S}(E)$, $\mathcal{R}(E)$ and $\mathcal{F}(E)$ respectively. Similarly, let $\mathcal{L}_I$, $\mathcal{J}_I$ and $\mathcal{I}_I$ denote the defining ideals of $\mathcal{S}(I)$, $\mathcal{R}(I)$ and $\mathcal{F}(I)$ respectively.
 
    By construction it then follows that $\,\displaystyle{\mathcal{L}_I = \mathcal{L}_E R'' + (X_1, \ldots, X_{e-1})},\,$ as well as $\,\displaystyle{\mathcal{J}_I = \mathcal{J}_E R'' + (X_1, \ldots, X_{e-1})}$. Since $S$ is standard graded, this implies that
      $$\mathcal{I}_I = [\mathcal{J}_I]_{(0,-)}= [\mathcal{J}_E]_{(0,-)} R'' + (X_1, \ldots, X_{e-1}) = \mathcal{I}_E R''+ (X_1, \ldots, X_{e-1}). $$
    Hence, $I$ is of fiber type if and only if 
      $$\mathcal{J}_I = \mathcal{L}_I + [\mathcal{J}_I]_{(0,-)} R'' =  \mathcal{L}_E + [\mathcal{J}_E]_{(0,-)} R'' + (X_1, \ldots, X_{e-1}).$$
    
    Now, if either assumption (i) or (ii) are satisfied, by \cref{backward} or \cref{NewSUV3.7} it follows that $\,X_1, \ldots, X_{e-1}\,$  form a regular sequence modulo $\mathcal{J}_E R''$. Therefore, it follows that 
      \begin{eqnarray*}
        \mathcal{J}_E R'' & = & \Big( \mathcal{L}_E R'' + [\mathcal{J}_E]_{(0,-)} R''+ (X_1, \ldots, X_{e-1}) \Big) \cap \mathcal{J}_E R'' \\
                          & = & \mathcal{L}_E R''   + [\mathcal{J}_E]_{(0,-)} R''+ (X_1, \ldots, X_{e-1}) \,\mathcal{J}_E R''
      \end{eqnarray*}
    By the graded version of Nakayama's Lemma, this means that 
       $$ \mathcal{J}_E R'' = \mathcal{L}_E R''   + [\mathcal{J}_E]_{(0,-)} R'', $$
    which can occur if and only if in $\,R[T_1, \ldots, T_n]\,$ one has 
      $\,\displaystyle{\mathcal{J}_E = \mathcal{L}_E   + [\mathcal{J}_E]_{(0,-)}}, \,$
    i.e. if and only if $E$ is of fiber type. $\blacksquare$

\subsection{Almost linearly presented modules of projective dimension one} 
 
 In this section we describe the Rees algebra and fiber cone of almost linearly presented modules of projective dimension one. Throughout we will consider the situation of \cref{SetDefEqs} below. 

 \begin{set} \label{SetDefEqs} \hypertarget{SetDefEqs}{}
    \em{Let $R=k[Y_1, \ldots, Y_d]$ be a polynomial ring over a field $k$, where $d \geq 2$. Let $E$ be a finite $R$-module, minimally generated by homogeneous elements of the same degree. Assume also that $E$ has projective dimension one and satisfies $G_d$. Then, has positive rank $e$ and admits a minimal free resolution of the form $\,\displaystyle{0 \to R^{\,n-e} \stackrel{\varphi}{\longrightarrow} R^n \to E \to 0},\,$ where $n= \mu(E)$. Assume that $\varphi$ is \emph{almost linear}, i.e. has linear entries, except possibly for those in the last column, which are homogeneous of degree $m \geq 1$.} 
 \end{set}
 
 In the situation of \cref{SetDefEqs}, after localizing at the unique homogeneous maximal ideal, by \cref{truetdefBourbaki} $\,E$ admits a generic Bourbaki ideal $I$, which is perfect of grade 2. Let $\psi$ be a minimal presentation of $I$ obtained from $\varphi$ as in \cref{BourbakiPres}. By construction, $\psi$ is also almost linear. In particular, the defining ideal of $\,\mathcal{R}(I)\,$ is described by the following theorem of Boswell and Mukundan \cite[5.3 and 5.6]{BM}. 

 \begin{thm} \label{BM} \hypertarget{BM}{}
    Let $R=k[Y_1, \ldots,Y_d]$ be a standard graded polynomial ring over a field $k$. Let $I$ be a perfect ideal of height 2 admitting an almost linear presentation $\psi$. Assume that $I$ satisfies $G_d$ and that $\mu(I)=d+1$. Then, the defining ideal of the Rees algebra $\mathcal{R}(I)$ is 
       $$\mathcal{J}= (\underline{Y} \cdot B(\psi)) + I_d(B_m(\psi)) = (\underline{Y} \cdot B(\psi)) \colon (Y_1, \ldots, Y_d)^m,$$
    where $m$ is the degree of the non-linear column of $\psi$ and $B_m(\psi)$ is the $m$th-iterated Jacobian dual of $\psi$ as in Definition~\ref{IterJacDuals}. Moreover:
    \begin{itemize}
       \item[$($i$)$] $\mathcal{R}(I)$ is almost Cohen-Macaulay, i.e. $\mathrm{depth}_{\,}\mathcal{R}(I) \geq d-1$, and it is Cohen-Macaulay if and only if $m=1$. 
       \item[$($ii$)$] $\,\mathcal{F}(I)$ is Cohen-Macaulay.
    \end{itemize}
 \end{thm}

 We now generalize  \cref{BM} to almost linearly presented modules of projective dimension one. 

 \begin{thm} \label{GenBM}
    Under the assumptions of \cref{SetDefEqs}, set $\underline{Y}= [Y_1, \ldots, Y_d]$ and assume that $\,n=d+e$. Then, the defining ideal of $\, \mathcal{R}(E)$ is 
       $$ \mathcal{J}= ((\underline{Y} \cdot B(\varphi))\, \colon (\underline{Y})^m) = (\underline{Y} \cdot B(\varphi)) + I_d(B_m(\varphi)),$$
    where $m$ is the degree of the non-linear column of $\varphi$ and $B_m(\varphi)$ denotes an $m$-th iterated Jacobian dual as in Definition~\ref{IterJacDuals}. Moreover:
    \begin{itemize}
       \item[$($i$)$] $\mathcal{R}(E)$ is almost Cohen-Macaulay, and it is Cohen-Macaulay if and only if $m=1$. 
       \item[$($ii$)$] $\,\mathcal{F}(E)$ is Cohen-Macaulay.
    \end{itemize}
 \end{thm}

    \emph{Proof}. We modify the proof of \cite[4.11]{SUV2003}. Let $a_1, \ldots, a_n$ be a minimal generating set for $E$ corresponding to the presentation $\,\varphi$, and let $ \,R[T_1, \ldots, T_n] \twoheadrightarrow \mathcal{R}(E) \,$ be the natural epimorphism, mapping $T_i$ to $a_i$ for all $i$. Localizing at the unique homogeneous maximal ideal, we may assume that $R$ is local and that $E$ admits a generic Bourbaki ideal $I$, which is perfect of grade 2 and such that $\, \mu(I)=n-e+1=d+1.\,$ If $e=1$, then $E \cong I$ and the statement follows from \cref{BM}. 

    So, assume that $e \geq 2$. With $x_j$ as in \cref{trueNotationBourbaki}, for $1 \leq j \leq e-1\,$ set $\,X_j= \sum_{i=1}^{n} Z_{ij} T_i, \,$ and note that $X_j$ is mapped to $x_j$ under the epimorphism $\,R''[T_1, \ldots, T_n] \twoheadrightarrow \mathcal{R}(E''). \,$ Set $\,\underline{T}= [T_1, \ldots, T_n]$. As in \cref{BourbakiPres}, we can construct a minimal almost linear presentation $\psi$ of $I$, such that 
      $$ [\underline{Y}] \cdot B(\varphi) \equiv [\underline{T}] \cdot  \left[ \begin{array}{c} 
         0 \\
         \hline
         \psi \\
      \end{array} \right]\, \; \mathrm{modulo} \, (X_1, \ldots, X_{e-1}). $$
    Let $B(\psi)$ be a Jacobian dual of $\psi$ defined by $\, \displaystyle  [\underline{T}] \cdot \left[ \begin{array}{c}
         0 \\
         \hline
         \psi \\
      \end{array} \right]  = [\underline{Y}] \cdot B(\psi)$. 
    Then, by \cref{BM}, the defining ideal of $\,\mathcal{R}(I)$ is 
     \begin{equation} \label{eqJI} \hypertarget{eqJI}{}
        \mathcal{J}_I = (\underline{Y} \cdot B(\psi))+ I_d(B_m(\psi)) = (\underline{Y} \cdot B(\psi)) \, \colon (\underline{Y})^m,
     \end{equation}
    where $m$ is the degree of the non-linear column of $\varphi$. Moreover, $\,\mathcal{R}(I)$ is almost Cohen-Macaulay, and Cohen-Macaulay if and only if the entries of $\psi$ are all linear, while $\,\mathcal{F}(I)$ is Cohen-Macaulay. In particular, by \cref{trueMainBourbaki} it follows that $\mathcal{R}(E)$ is Cohen-Macaulay if and only if $\varphi$ is linear. 

    To prove the remaining statements, notice that $E''_{\mathfrak{q}}\,$ is of linear type for all primes $\mathfrak{q}$ in the punctured spectrum of $R''$ (this is because $E''$ has projective dimension one and satisfies $G_d$, by \cite[Propositions 3 and 4]{Avramov}). Hence, also $I_{\mathfrak{q}}\,$ is of linear type for the same primes $\mathfrak{q}$. Moreover, the discussion above shows that $\, \displaystyle{ \mathrm{depth}_{\,}\mathcal{R}(I) \geq \mathrm{dim}_{\,}\mathcal{R}(I) -1=d \geq 2}$. Hence, inducting on $e\,$ and using \cref{NewSUV3.7}$\,$ in the case when $e=2$, we obtain that $\, \displaystyle{\mathcal{R}(I) \cong \mathcal{R}(E'')/(F'')} \,$ and $\,x_1, \ldots, x_{e-1}\,$ form a regular sequence on $\,\mathcal{R}(E'')$. Thus, $X_1, \ldots, X_{e-1}$ form a regular sequence modulo $\mathcal{J} R''$. This shows that $\,\mathcal{R}(E'')$ is almost Cohen-Macaulay, whence $\mathcal{R}(E)$ is almost Cohen-Macaulay. It also implies that 
       $$ \mathcal{J}_I= \mathcal{J} R'' + (X_1, \ldots, X_{e-1}) $$
    Hence, from (\ref{eqJI}) and \cref{IterJacDualsPass} below it follows that
       $$ \mathcal{J}R'' + (X_1, \ldots, X_{e-1})= (\underline{Y} \cdot B(\varphi)) + I_d(B_m(\varphi)) +(X_1, \ldots, X_{e-1}). $$

   On the other hand, since $E$ is of linear type locally on the punctured spectrum of $R$, it follows that
     $$ \mathcal{J} \supseteq (\underline{Y} \cdot B(\varphi)) \, \colon (\underline{Y})^m \supseteq (\underline{Y} \cdot B(\varphi)) + I_d(B_m(\varphi)), $$
   where the last inclusion follows from Theorem~\ref{IterJacDuals}(c). Therefore, since $X_1, \ldots, X_{e-1}$ form a regular sequence modulo $\mathcal{J}R'',\,$ in $R''[T_1, \ldots, T_n]$ we have:
    \begin{eqnarray*}
      \mathcal{J} \!\!\! & = & \!((\underline{Y} \cdot B(\varphi)) + I_d(B_m(\varphi)) +(X_1, \ldots, X_{e-1})) \cap \mathcal{J} \\
      \, & =  & ((\underline{Y} \cdot B(\varphi)) + I_d(B_m(\varphi))) +  (X_1, \ldots, X_{e-1})\, \mathcal{J} .
    \end{eqnarray*}
   By the graded version of Nakayama's Lemma, this means that 
     $$\mathcal{J} = (\underline{Y} \cdot B(\varphi)) + I_d(B_m(\varphi))= (\underline{Y} \cdot B(\varphi)) \, \colon (\underline{Y})^m, $$
   as claimed. Finally, since $\mathcal{F}(I)$ is Cohen-Macaulay and $\,  \mathrm{depth}_{\,}\mathcal{R}(I) \geq 2$, from \cref{3.5FiberCone}(b) it follows that $\,\mathcal{F}(E)\,$ is Cohen-Macaulay.   $\blacksquare$\\

 \begin{lemma} \label{IterJacDualsPass} \hypertarget{IterJacDualsPass}{}
   Let $R=k[Y_1, \ldots, Y_d]_{(Y_1, \ldots, Y_d)}$, and denote $\underline{Y}= [Y_1, \ldots, Y_d]$. Let $\varphi$, $\psi$, $B(\psi)$, and $X_1, \ldots, X_{e-1}\,$ be as in the proof of \cref{GenBM}. Then, for all $i$ and for any Jacobian dual $\,B(\varphi)$ of $\varphi$, in $\,R''[T_1, \ldots, T_n]$ we have 
     $$(\underline{Y} \cdot B(\varphi)) + I_d(B_i(\varphi)) +(X_1, \ldots, X_{e-1})= (\underline{Y} \cdot B(\psi))+ I_d(B_i(\psi)).$$
 \end{lemma}

   \emph{Proof}. Choose $B(\psi)$ such that $[\underline{Y}] \cdot B(\psi)= [\underline{T}] \cdot  \left[ \begin{array}{c}
     0 \\
     \hline
     \psi \\
     \end{array} \right], \,$ as in the proof of \cref{GenBM}. Then, in $R''[T_1, \ldots, T_n]$ we have
     $$ [\underline{Y}] \cdot B(\varphi) \equiv [\underline{Y}] \cdot B(\psi) \: \mathrm{modulo} \,(X_1, \ldots, X_{e-1}). $$
   So, the statement is proved for $i=1$.  Now, let $i+1 \geq 2$ and assume that the statement holds for $B_i(\varphi)$. Let $C_i$ be a matrix as in Definition~\ref{IterJacDuals}. Since 
     $$(\underline{Y} \cdot B_i(\varphi)) + (I_d(B_i(\varphi)) \cap (\underline{Y})) = (\underline{Y} \cdot B_i(\varphi)) + (\underline{Y} \cdot C_i)$$
   and the $B_i(\varphi)$ are bigraded, going modulo $(X_1, \ldots, X_{e-1})$, in $R''[T_1, \ldots, T_n]$ we have 
     $$(\underline{Y} \cdot B_i(\psi)) + (I_d(B_i(\psi)) \cap (\underline{Y})) = (\underline{Y} \cdot B_i(\psi)) + (\underline{Y} \cdot \overline{C_i}),$$
   where $\overline{C_i}$ denotes the image of $C_i$ modulo $(X_1, \ldots, X_{e-1})$. Now, let $\,B_{i+1}(\psi) = [B_i(\psi)) \, | \, \overline{C_i}]$. Then, in $R''[T_1, \ldots, T_n]\,$ we have that
     $$(\underline{Y} \cdot B(\varphi)) + I_d(B_{i+1}(\varphi)) +(X_1, \ldots, X_{e-1})= (\underline{Y} \cdot B(\psi))+ I_d(B_{i+1}(\psi)),$$
  as we aimed to show. $\, \blacksquare$\\

  We remark that the equality $\, \mathcal{J}= (\underline{Y} \cdot B(\varphi)) \colon (Y_1, \ldots, Y_d)^m \, $ for the defining ideal of the Rees algebra of a module $E$ as in \cref{GenBM} could be obtained without using generic Bourbaki ideals, by modifying the proof of \cite[6.1(a)]{KPU}. In fact, even if \cite[6.1(a)]{KPU} is stated for perfect ideals of height two, its proof only uses the structure of the presentation matrix, and one would only need to  adjust the ranks to prove the statement for modules of projective dimension one. More generally, up to this minor adjustment in the proof, \cite[6.1(a)]{KPU} shows that if $E$ is a finite module over $\,k[Y_1, \ldots, Y_d]\,$ minimally generated by homogeneous elements of the same degree, then the defining ideal of $\mathcal{R}(E)$ is $\, \mathcal{J}= (\underline{Y} \cdot B(\varphi)) \, \colon (\underline{Y})^N, \, $ where $\, N= 1+ \sum_{i=1}^d (\epsilon_i -1) \,$ and the $\epsilon_i$ are the degrees of the columns of $\varphi$. 
  
  Similarly, a good portion of the proof of \cref{BM} could be adjusted to the case of modules of projective dimension one by modifying the ranks of the presentation. However, we would not be able to generalize the whole statement of \cref{BM} using this method. Indeed, the proof of the equality $\, \mathcal{J}= (\underline{Y} \cdot B(\varphi)) + I_d(B_i(\varphi))\,$ in the case of perfect ideals of height two crucially makes use of the ideal structure of the cokernel of the presentation matrix (see the proof of \cite[5.3]{BM}).

\section*{Acknowledgements} 
Most of the work presented in this manuscript was part of the author's Ph.D. thesis. The author wishes to thank her advisor Bernd Ulrich for his insightful comments on some of the results here presented and for assigning \cref{homework} as a homework problem in one of his courses. Also, part of the content of \cref{SecFiberCone} was motivated by a question of Jonathan Monta\~no, whom we thank very much for fruitful conversations on the topic. 

\end{document}